\documentclass{amsart}[12pt]

\usepackage{amsmath,amsfonts,amssymb,amsthm}
\usepackage{enumitem}
\usepackage[explicit]{titlesec}

\makeatletter
\def\@secnumpunct{. }
\makeatother

\titleformat{\subsection}[block]
{}{\thetitle. }{0pt}{\bf{#1.}}

\titleclass{\subsectionn}{straight}[\section]
\newcounter{subsectionn}

\titleformat{\subsectionn}[runin]{}{\thetitle. }{0pt}{}
\titlespacing{\subsectionn}{0pt}{*1.5}{*1.5}

\newtheorem{thm}{Theorem}[section]
\newtheorem{lemma}[thm]{Lemma}

\newtheorem{cor}[thm]{Corollary}

\newtheorem{conj}[thm]{Conjecture}

\theoremstyle{definition}

\def\zz{\mathbb Z}

\def\gg{\mathbb G}

\def\wt{\widetilde}
\def\ol{\overline}
\def\<{\langle}
\def\>{\rangle}

\def\0{{\mathbf 0}}

\newcommand{\ra}{\rightarrow}

\newcommand{\lb}{\langle}
\newcommand{\rb}{\rangle}

\newcommand{\deq}{\mathrel{\mathop:}=}

\newcommand{\pathto}{\:\rightsquigarrow\:}
\newcommand{\edgeto}{\leftrightarrow}
\newcommand{\SL}{SL}

\newcommand{\setc}[3][]{\left\{ #2 \,\middle|\, #3 \vphantom{#1|}\right\}}

\newcommand{\gen}[1]{\left\lb #1 \right\rb}
\newcommand{\abs}[2][]{\left|#2 \vphantom{#1|}\right|}

\DeclareMathOperator{\Aut}{Aut}
\DeclareMathOperator{\Out}{Out}

\DeclareMathOperator{\Cay}{Cay}

\DeclareMathOperator{\im}{image}

\begin{document}

\title{Non-amenabilty of product replacement graphs}

\author[Anton~Malyshev]{ \ Anton~Malyshev$^\star$}

\thanks{\thinspace ${\hspace{-.45ex}}^\star$Department of Mathematics,
UCLA, Los Angeles, CA, 90095.
\hskip.06cm
Email:
\hskip.06cm
\texttt{amalyshev@math.ucla.edu}}

\begin{abstract}
We prove non-amenability of the product replacement graphs $\Gamma_n(G)$ for uniformly non-amenable groups. We also prove it for $\zz$-large groups, when $n$ is sufficiently large. It follows that $\Gamma_n(G)$ is non-amenable when $n$ is sufficiently large for hyperbolic groups, linear groups, and elementary amenable groups.
\end{abstract}

\maketitle
\theoremstyle{plain}

\section{Introduction}

The \emph{product replacement graph} of a group~$G$ is the graph of generating $n$-tuples of~$G$, connected by Nielsen moves. These graphs play an important role in computational group theory (see e.g.~\cite{BL}), and have been considered in connection with the Andrews-Curtis conjecture (see e.g.~\cite{BLM,Met}). Relatively little is known about these graphs. Even their connectivity is a major open problem (see e.g.~\cite{evans,myropolska}). In this paper, we continue the investigation in~\cite{expgrowth}.

The main subject of this paper is the non-amenability of product replacement graphs. This is related to the well-known problem of whether the automorphism group $\Aut(F_n)$ of the free group has Kazhdan property $(T)$ for~$n > 3$, (see~\cite{property_tau,LP}, see also Subsection~\ref{subsec:fr:gamma_schreier}). If so, product replacement graphs of infinite groups would be \emph{non-amenable}. This motivates the following conjecture.

\begin{conj}\cite{expgrowth}\label{conj:all_nonamenable}
The product replacement graph $\Gamma_n(G)$ of an infinite finitely generated group $G$ is non-amenable for sufficiently large~$n$.
\end{conj}

The purpose of this note is to prove that Conjecture~\ref{conj:all_nonamenable} holds for several classes of groups. We prove and use the following two theorems.

\begin{thm}\label{thm:uniformly_nonamenable}
If $G$ is uniformly non-amenable, then the product replacement graph $\Gamma_n(G)$ is non-amenable for every~$n \geq d(G)$.
\end{thm}

We say $G$ is \emph{$\zz$-large} if $G$ has a finite index subgroup which has $\zz$ as a quotient.

\begin{thm}\label{thm:zz_large}
If $G$ is $\zz$-large, then the product replacement graph $\Gamma_n(G)$ is non-amenable for sufficiently large~$n$.
\end{thm}

We combine these theorems with the results on uniform non-amenability in~\cite{uniform_nonamenability} to show that several classes of infinite groups have non-amenable product replacement graphs. In particular, hyperbolic groups, linear groups, elementary amenable groups, and free Burnside groups all satisfy Conjecture~\ref{conj:all_nonamenable}. Note that elementary amenable groups include virtually solvable groups and virtually amenable groups.

\smallskip

The paper is structured as follows. In Section~\ref{section:definitions}, we define our terms and recall some basic facts about non-amenability. In Sections~\ref{section:zz_large}~and~\ref{section:uniformly_nonamenable}, we prove our main theorems. In Section~\ref{section:classes} we discuss the consequences of these theorems to several general classes of groups. In Section~\ref{section:lemma_proofs}, we prove the lemmas we used in the previous sections. Finally, in Section~\ref{section:final_remarks} we discuss the relationship of this problem to unsolved problems, and indicate further directions.

\section{Definitions}\label{section:definitions}

The \emph{product replacement graph} $\Gamma_n(G)$ of a finitely generated group $G$ is the undirected graph whose vertices are $n$-tuples $S = (s_1, \dots, s_n) \in G^n$ with $G = \gen{s_1, \dots, s_n}$, and edges
\begin{align*}
(s_1, \dots, s_i, \dots, s_j, \dots s_n) & \edgeto (s_1, \dots, s_i, \dots, s_js_i^{\pm 1}, \dots, s_n)
\\
(s_1, \dots, s_i, \dots, s_j, \dots s_n) & \edgeto (s_1, \dots, s_i, \dots, s_i^{\pm 1} s_j, \dots, s_n)
\end{align*}
for each pair~$1 \leq i, j \leq n$, with~$i \neq j$. We call a step from a vertex in this graph to one of its neighbors a \emph{Nielsen move}. We denote the minimum number of generators of $G$ with~$d(G)$, so $\Gamma_n(G)$ is nonempty only for~$n \geq d(G)$.

Here and elsewhere, by a \emph{graph} $\Gamma = (V,E)$ we mean a possibly infinite undirected graph of bounded degree, which may have loops and repeated edges. By an abuse of notation, we often write $\Gamma$ to mean the vertex set of~$\Gamma$. The \emph{distance} $d(v,w)$ between two vertices $v,w \in \Gamma$ is defined to be the length of the shortest path connecting them (or $\infty$ if there is no such path). We define $d(v,X)$ to be the distance from a vertex $v \in \Gamma$ to a set of vertices~$X \subset \Gamma$, i.e.~$d(v,X) = \min \{d(v,x) \mid x \in X\}$. For any two sets of vertices $X,Y \subset \Gamma$, we write $E(X,Y)$ for the set of edges between $X$ and $Y$.

Let $G$ be a finitely generated group, and fix a generating tuple~$S = (s_1, \dots, s_n) \in \Gamma_n(G)$. The \emph{Cayley graph} ${\Cay(G,S) = (V,E)}$ is the graph with vertex set~$V = G$, and edges $g \edgeto g s_i^{\pm 1}$ for each~$1 \leq i \leq n$.

We define the \emph{Cheeger constant} of a nonempty graph $\Gamma = (V,E)$ by
$$h(\Gamma)
\deq \inf_X \frac{\abs{\partial X}}{\abs{X}},$$
where the infimum runs over all finite sets of vertices~$X \subset \Gamma$, and $\partial X = E(X, \ol X)$ denotes the set of edges leaving~$X$. We say $\Gamma$ is \emph{non-amenable}\footnote{Many authors require a non-amenable graph to be connected, but we omit this requirement. Our definition still forbids finite connected components in~$\Gamma$.} if~$h(\Gamma) > 0$. When $\Gamma = \Cay(G,S)$, we denote the Cheeger constant $h(\Cay(G,S))$ with~$h(G,S)$. It is easy to check that the property $h(G,S) > 0$ depends only on the group $G$ and not on the choice of generators~$S$. We say a finitely generated group $G$ is \emph{non-amenable} if $h(G,S)>0$, i.e.\ if $\Cay(G,S)$ is non-amenable. However, the Cheeger constant $h(G,S)$ itself may depend on~$S$. We say~$G$ is \emph{uniformly non-amenable}\footnote{This is weaker than the definition in~\cite{uniform_nonamenability}, as we do not compare non-amenability with respect to generating tuples of different lengths.} if for every $n \geq d(G)$, there is a constant $h(G,n) > 0$ such that $$\inf_{S \in \Gamma_n(G)} h(G,S) > h(G,n).$$

We say a map $f:\Gamma' \ra \Gamma$ between two graphs is a \emph{Lipschitz} map if there is a constant $C>0$ such that for every pair of neighbors $v, w \in \Gamma'$, we have $d(f(v),f(w)) \leq C$. We say a subset $W \subset \Gamma$ is \emph{dense} if there is a constant $D>0$ such that for every $v \in \Gamma$ we have $d(v,W) \leq D$.

\section{Uniformly non-amenable groups}\label{section:uniformly_nonamenable}

A natural special case of Conjecture~\ref{conj:all_nonamenable} is the case where $G$ is assumed to be a non-amenable group. One obstacle to proving this is that the Cheeger constant of the Cayley graph may be arbitrarily small, depending on the generating set of $G$ (see e.g.~\cite{uniform_nonamenability,weakly_amenable}). Theorem~\ref{thm:uniformly_nonamenable} asserts that if this is not the case, i.e.\ if the group is uniformly non-amenable, then $\Gamma_n(G)$ is non-amenable.

\subsection{Lemmas}

To prove this theorem, we need the following two lemmas. The first lemma is a variation on the fact that \emph{quasi-isometry} of graphs preserves non-amenability. For more on quasi-isometry, see e.g.~\cite[$\S$3,4]{woess},~\cite[$\S$IV.B]{de_la_harpe}.
\begin{lemma}\label{lem:transfer_lip}
Let $\Gamma$ and $\Gamma'$ be infinite graphs, with $\Gamma$ non-amenable. Let $f:\Gamma \ra \Gamma'$ be an injective Lipschitz map such that $f(\Gamma)$ is dense in~$\Gamma'$. Then $\Gamma'$ is also non-amenable.
\end{lemma}

The second lemma is related to the following standard fact: The class of amenable groups is closed under the operation of group extensions. That is, if $G$ is non-amenable, then for every normal subgroup $H$ of~$G$, either $H$ or~$G/H$ is non-amenable. In our argument, we only consider the case where $H$ is the center of $G$, but we need an explicit lower bound on the Cheeger constant of~$G/H$:
\begin{lemma}\label{lem:transfer_centermod}
If $G$ is a non-amenable group, then $G/Z(G)$ is also non-amenable. Moreover, for every $S = (s_1, \dots, s_n) \in \Gamma_n(G)$, we have $$h(G/Z(G),\wt S) \geq h(G,S),$$ where $\wt S = (s_1Z(G), \dots, s_nZ(G)) \in \Gamma_n(G/Z(G))$.
\end{lemma}

We prove these lemmas in Section~\ref{section:lemma_proofs}.

\subsection{Proof of Theorem~\ref{thm:uniformly_nonamenable}}

Let $G$ be a uniformly non-amenable group. Given any $S = (s_1, \dots, s_n) \in \Gamma_n(G)$, define a map $f_S:G \ra \Gamma_n(G)$ by
$$
f_S(g) = gSg^{-1} = (gs_1g^{-1}, \dots, gs_ng^{-1})
$$
Observe that
\begin{align*}
f_S(g s_1)
&
= \Big(gs_1 s_1s_1^{-1} g^{-1},\: g s_1 s_2 s_1^{-1} g^{-1},\: \dots \:,\: g s_1 s_n s_1^{-1} g^{-1}\Big)
\\&
= \Big(g s_1 g^{-1},\: (g s_1 g^{-1}) (g s_2 g^{-1}) (g s_1 g^{-1})^{-1},\: \dots\: ,\:  (g s_1 g^{-1}) (g s_n g^{-1}) (g s_1 g^{-1})^{-1} \Big)
\end{align*}
is within $2n-2$ Nielsen moves of~$f_S(g)$. The same is true of $f_S(g s_i)$ for every $1 \leq i \leq n$, so $f_S$ is a Lipschitz map from $\Cay(G,S)$ to~$\Gamma_n(G)$, with Lipschitz constant~$2n-2$.

Let $\wt G = G/Z(G)$, where $Z(G)$ denotes the center of $G$. For each $g \in G$, denote by $\wt g$ the projection of~$g$ into~$G/Z(G)$. Let $\wt S$ denote the corresponding generating $n$-tuple of~$\wt G$, i.e.~$\wt S = (\wt s_1, \dots, \wt s_n)$. Observe that $f_S(g) = f_S(h)$ if and only if $g s_i g^{-1} = h s_i h^{-1}$ for every $1 \leq i \leq n$. This occurs precisely when $g^{-1}h$ commutes with every~$s_i$, i.e.\ when $hg^{-1} \in Z(G)$. Thus we have a well-defined and injective induced map $\wt f_S: \wt G \ra \Gamma_n(G)$ given by $\wt f_S(\wt g) = f_S(g)$. We also have that $\wt f_S$ is a Lipschitz map from $\Cay(\wt G, \wt S)$ to $\Gamma_n(G)$, with Lipschitz constant~$2n-2$.

Now let $S$ vary. Given $S' \in \Gamma_n(G)$, note that $S' \in \im \wt f_S$ if and only if $S' = gSg^{-1}$ for some~$g \in G$. This is an equivalence relation, so the images of the maps $\wt f_S$ form a partition of $\Gamma_n(G)$ into equivalence classes. Let $\mathcal S$ be a set of representatives of these equivalence classes, and consider the disjoint union of graphs
$$\Delta \deq \coprod_{S \in \mathcal S} \Cay(\wt G, \wt S).$$
Using Lemma~\ref{lem:transfer_centermod} and the fact that $G$ is uniformly non-amenable, we can conclude that the graph $\Delta$ is non-amenable, since
$$h(\Delta) \geq \inf_{S \in \mathcal S} h(\wt G, \wt S) \geq \inf_{S \in \mathcal S} h(G, S) > 0$$
The maps $\wt f_S$ for $S \in \mathcal S$ combine into one map
$\wt f: \Delta \ra \Gamma_n(G).$
This map is injective, surjective, and Lipschitz. Therefore, by Lemma~\ref{lem:transfer_lip}, $\Gamma_n(G)$ is also non-amenable.
\qed

\section{$\zz$-large groups}\label{section:zz_large}

It was shown in \cite{expgrowth} that a group containing an element of infinite order must have exponentially growing product replacement graphs $\Gamma_n(G)$, for sufficiently large $n$. In order to guarantee non-amenability, we require a stronger property: following~\cite{largeness}, we say a group $G$ is \emph{$\zz$-large} if it contains a finite index subgroup which has $\zz$ as a quotient. Theorem~\ref{thm:zz_large} asserts that $\zz$-large groups have non-amenable product replacement graphs $\Gamma_n(G)$, for sufficiently large~$n$.

Note that in this case there may exist $n \geq d(G)$ for which $\Gamma_n(G)$ fails to be non-amenable. For example, the infinite dihedral group $D_\infty$ is $\zz$-large, but $\Gamma_2(D_\infty)$ is an amenable infinite graph. However, $\Gamma_n(D_\infty)$ is non-amenable for every~$n \geq 3$.

\subsection{Lemmas}

To prove Theorem~\ref{thm:zz_large}, we again use Lemma~\ref{lem:transfer_lip}, as well as two additional lemmas. If $H$ is a quotient of~$G$, then $\Gamma_n(G)$ is a lift of~$\Gamma_n(H)$. A lift of a non-amenable graph is non-amenable, so we have:
\begin{lemma}\label{lem:transfer_lift}
Let $G$ be a finitely generated group. If $\Gamma_n(H)$ is non-amenable for some quotient $H$ of $G$, then $\Gamma_n(G)$ is non-amenable.
\end{lemma}

The following fact is well-known. It is related to the fact that $\Gamma_2(\zz^k)$ is a Schreier graph of $\SL(2,\zz)$, which has property $(\tau)$ with respect to its congruence subgroups (see~\cite{property_tau}).
\begin{lemma}\label{lem:zz}
The product replacement graph $\Gamma_2(\zz^k)$ is non-amenable for every~$k > 0$.
\end{lemma}

We refer the reader to Section~\ref{section:lemma_proofs} for the proofs.

\subsection{Proof of Theorem~\ref{thm:zz_large}}

Let $G$ be a finitely generated $\zz$-large group. Let $H$ be a finite index subgroup of $G$ which has $\zz$ as a quotient. Then $H$ is also finitely generated. Thus, the statement that $H$ has $\zz$ as a quotient is equivalent to the statement that $[H:H'] = \infty$, where $H' = [H,H]$ is the commutator subgroup of~$H$. This, in turn, is equivalent to the statement that $[G:H'] = \infty$. Let $H^\circ = \bigcap_{g \in G} gHg^{-1}$ be the normal core of $H$ in $G$. Then we also have $[G: H^\circ] < \infty$ and $[G:(H^\circ)'] = \infty$. That is, $H^\circ$ satisfies the same hypotheses as $H$, so by replacing $H$ with $H^\circ$ if necessary, we may assume that $H$ is normal in~$G$.

We have that $H/H' \cong \zz^r \times A$ for some positive integer $r$ and some finite abelian group~$A$. Let~$N$ be the kernel of the corresponding homomorphism $H \ra \zz^r$. Then $N$ is a characteristic subgroup of~$H$, and therefore $N$ is a normal subgroup of~$G$.

By Lemma~\ref{lem:transfer_lift}, it is enough to show that $\Gamma_n(G/N)$ is non-amenable for sufficiently large~$n$. Thus, by replacing $G$ with $G/N$ and $H$ with $H/N$, we may assume $N$ is trivial, and hence $H \cong \zz^r$.

Fix $n \geq \log_2 \abs{G/H}$. Consider a generating tuple $S = (s_1, \dots, s_{n+2}) \in \Gamma_{n+2}(G)$, and the corresponding tuple $\wt S = (\wt s_1, \dots, \wt s_{n+2}) \in \Gamma_n(G/H)$. Every generating $(n+1)$-tuple of $G/H$ is redundant (see e.g.~\cite[2.2]{whatprp}), so a bounded number of Nielsen moves in $\Gamma_n(G/H)$ sends
$$(\wt s_1, \dots, \wt s_{n+2}) \pathto (\wt t_1, \dots, \wt t_n, 1, 1).$$ The same Nielsen moves in $\Gamma_n(G)$, then, send
$$(s_1, \dots, s_{n+2}) \pathto (t_1, \dots, t_n, h_1, h_2),$$
where $h_1, h_2 \in H$. If $h_1$ and $h_2$ are both trivial, then $G = \gen{t_1, \dots, t_n}$, so in at most $[G:H]$ more Nielsen moves, we can reach an element of that form with $h_1$, $h_2$ not both trivial.

For every nontrivial subgroup $K < H$ and every $T = (t_1, \dots, t_n) \in G^n$ with $G = \gen{T,K}$, there is a graph embedding $\Gamma_2(K) \ra \Gamma_{n+2}(G)$ given by $$(h_1, h_2) \mapsto (t_1, \dots, t_n, h_1, h_2),$$ and the images of these embedding are disjoint. Let $\Delta$ denote the union of these embeddings. In the previous paragparh, we showed that every vertex of $\Gamma_{n+2}(G)$ is a bounded distance away from $\Delta$. Since each $K$ satisfies $K \cong \zz^k$ for some $1 \leq k \leq r$, we have
$$h(\Delta) \geq \inf_{T,K} h(\Gamma_2(K)) \geq \min_{1 \leq k \leq r} h\left(\Gamma_2(\zz^k)\right) > 0.$$
Thus, $\Delta$ is non-amenable. Lemma~\ref{lem:transfer_lip} implies that $\Gamma_{n+2}(G)$ is also non-amenable.
\qed

\section{Examples}\label{section:classes}

Theorems~\ref{thm:uniformly_nonamenable} and \ref{thm:zz_large} combine to show that several nice classes of groups satisfy Conjecture~\ref{conj:all_nonamenable}. First of all, there are a number of uniformly non-amenable groups, as shown in~\cite{uniform_nonamenability}:

\begin{thm}[\cite{uniform_nonamenability}]\label{thm:uniformly_nonamenable_classes}
The following classes of groups are uniformly non-amenable:
\begin{enumerate}[label=\emph{(\roman*)}]
\item \label{itm:first} non-elementary word-hyperbolic groups,
\item large groups (i.e.\ groups with a finite index subgroup which has $F_2$ as a quotient),
\item free Burnside groups $B(m,n)$ with $m \geq 2$ and sufficiently large odd $n$,
\item \label{itm:last} groups which act acylindrically on a simplicial tree without fixed points.
\end{enumerate}
\end{thm}

Combining this result with Theorem~\ref{thm:uniformly_nonamenable}, we obtain:
\begin{cor}\label{cor:every_n}
The product replacement graph $\Gamma_n(G)$ is non-amenable for every $n \geq d(G)$, if $G$ belongs to one of the classes \ref{itm:first}-\ref{itm:last} of Theorem~\ref{thm:uniformly_nonamenable_classes}.
\end{cor}

Using Theorem~\ref{thm:zz_large}, we can extend this result to a larger class of groups, at the cost of a somewhat weaker conclusion. First, we make the following observation, which we prove in Section~\ref{section:lemma_proofs}:
\begin{lemma}\label{lem:EA_zz-large}
Let $G$ be an infinite finitely generated group. If $G$ is elementary abelian, then $G$ is $\zz$-large.
\end{lemma}

By Gromov's Theorem~\cite{gromov}, every infinite group of polynomial growth is virtually nilpotent, and therefore elementary amenable. By definition, infinite elementary hyperbolic groups contain $\zz$ as a finite index subgroup. Finally, every virtually solvable group is elementary amenable. 

We also have the following theorem, which follows from a stronger version of the Tits alternative proven in~\cite{linear_independence}.
\begin{thm}[{\cite[Theorem 1.5]{linear_independence}}]\label{thm:linear_uniformly_nonamenable}
If $G$ is a linear group, then either $G$ is virtually solvable, or $G$ is uniformly non-amenable.
\end{thm}

Combining these observations with Corollary~\ref{cor:every_n}, we obtain:
\begin{cor}\label{cor:sufflarge_n}
The product replacement graph $\Gamma_n(G)$ is non-amenable for sufficiently large $n$, if $G$ is an infinite finitely generated group which belongs to one of the following classes:
\begin{enumerate}[label=\emph{(\roman*)}]
\item elementary amenable groups,
\item groups of polynomial growth,
\item word-hyperbolic groups,
\item linear groups.
\end{enumerate}
\end{cor}

\section{Proofs of Lemmas}\label{section:lemma_proofs}

We now prove the lemmas used in the previous sections. The arguments in this section are standard, but we need the results in a specific form.

\subsection{Proof of Lemma~\ref{lem:transfer_lip}}

Let $\Gamma$ and $\Gamma'$ be any infinite graphs, where $\Gamma$ is non-amenable. Let $f: \Gamma \ra \Gamma'$ be an injective Lipschitz map with Lipschitz constant~$C$. Suppose that $d(x,f(\Gamma)) \leq D$ for every $x \in \Gamma'$.

Given a finite set of vertices $X \subset \Gamma'$, define the \emph{$r$-neighborhood} of $X$ to be $$X^{(r)} = \{v \in \Gamma' \mid d(v,X) \leq r\}.$$ Let $d \geq 2$ be an upper bound on the degrees of vertices in $\Gamma'$ and $\Gamma$. Suppose $\abs{X^{(r)}} \geq \alpha \abs{X}$ for some~$\alpha > 1$. Then there are at least $(\alpha-1) \abs{X}$ paths of length $r$ or less from $X$ to $\ol X$, each of which contains at least one edge leaving~$X$. Each such edge occurs in at most $r d^{r-1} + (r-1) d^{r-2} + ... + 1 \leq r^2d^{r-1}$ of these paths, so $$\abs{\partial X} \geq \frac{\alpha - 1}{r^2 d^{r-1}} \abs{X}.$$
Thus, it is enough to show that there is positive integer $r$ and a constant $\alpha > 1$ such that $\abs{X^{(r)}} \geq \alpha \abs{X}$ for every finite subset~$X \subset \Gamma'$.

Let $C$ and $D$ be as above. Given a finite $X \subset \Gamma'$, every vertex of $X$ is within $D$ steps of some $v \in f(\Gamma)$, and for each $v \in f(\Gamma)$ there are at most $d^D + d^{D-1} + \dots + 1 \leq d^{D+1}$ vertices within $D$ steps of~$v$. It follows that
$$\abs{f^{-1}\big(X^{(D)}\big)} \geq \abs{X^{(D)} \cap f(\Gamma)} \geq \abs{X}/d^{D+1}.$$
Since $\Gamma$ is non-amenable, there are at least $h(\Gamma) \abs{X}/d^{D+1}$ edges leaving $f^{-1}(X^{(D)})$, and therefore at least $h(\Gamma) \abs{X}/d^{D+2}$ vertices $v \in \Gamma$ with $d(v,f^{-1}(X^{(D)})) = 1$. Each such $v$ maps to a unique $v' \in \Gamma'$ with $v' \notin X^{(D)}$ and $d(v', X^{(D)}) \leq C$. Hence,
$$\abs{X^{(D+C)}}
\geq \abs{X^{(D)}} + h(\Gamma) \abs{X}/d^{D+2}
\geq (1+h(\Gamma)/d^{D+2}) \abs{X},
$$
as desired.
\qed

\subsection{Proof of Lemma~\ref{lem:transfer_centermod}}

Let $G$ be a non-amenable group, and fix a generating $n$-tuple $S = (s_1, \dots, s_n)$ of~$G$. Let~$H = Z(G)$. Define $\wt G = G/H$, and $\wt S = (\wt s_1, \dots, \wt s_n) = (s_1H, \dots, s_nH)$. Let $\pi:G \ra \wt G$ denote the usual projection.

By picking representatives for each coset of $H$, we have a bijection between $G$ and $\wt G \times H$. That is, elements of $G$ can be represented in the form $(g, h) \in \wt G \times H$, where group operation is given by $(g_1, h_1)(g_2, h_2) = (g_1g_2, \star )$. In fact, because elements of $H$ commute with everything, we must have
$$(g_1, h_1)(g_2, h_2) = (g_1 g_2, \phi(g_1, g_2) h_1 h_2),$$ for some function $\phi:\wt G \times \wt G \ra H$. Then we can write the original generators $s_i$ in this form: $s_i = (\wt s_i, t_i)$, for some $t_i \in H$.

It is enough to show that for every finite subset $X \subset \Cay(\wt G, \wt S)$, we have $\abs{\partial X}/\abs{X} \geq h(G,S)$. Let $X$ be a finite subset of $\Cay(\wt G, \wt S)$, and let $Y = \pi^{-1}(X)$. Let $K$ be the subgroup of $H$ generated by $$T = \Big(\phi(x, \wt s_i)t_i\Big)_{\substack{x \in X \\ 1 \leq i \leq n.}}$$ Then $K$ is an abelian group, which implies it is amenable. Therefore, we have a sequence of finite subsets $B_1, B_2, \dots \subseteq \Cay(K, T)$ with $\abs{\partial B_k}/\abs{B_k} \ra 0$ as $k \ra \infty$. Let $C_k = \setc{(x,h)}{x \in X, h \in B_k}$.

Partition the set $\partial C_k$ of edges leaving $C_k$ into two sets: $$\partial_{\text{out}} C_k = E(C_k, G \setminus Y) \qquad \text{ and } \qquad \partial_{\text{in}} C_k = E(C_k, Y \setminus C_k).$$
We have that $\abs{\partial_{\text{out}}(C_k)} = \abs{B_k}\abs{\partial{X}}$ and $\abs{\partial_{\text{in}}(C_k)} \leq \abs{X}\abs{\partial B_k}$. Thus
\begin{align*}
h(G,S)
\leq
\frac{\abs{\partial C_k}}
{\abs{C_k}}
=
\frac{\abs{\partial_{\text{out}}(C_k)} + \abs{\partial_{\text{in}}(C_k)}}
{\abs{C_k}}
\leq
\frac{\abs{B_k}\abs{\partial X} + \abs{X}\abs{\partial B_k}}
{\abs{X}\abs{B_k}}
=
\frac{\abs{\partial X}}{\abs{X}} + \frac{\abs{\partial B_k}}
{\abs{B_k}}.
\end{align*}
Since the second term goes to $0$ as $k \ra \infty$, we have $h \leq \abs{\partial X}/{\abs{X}}$, as desired.

\subsection{Proof of Lemma~\ref{lem:transfer_lift}}
There is a characterization of non-amenability in terms of recurrent walks. Let $\Gamma = (V,E)$ be a nonempty $d$-regular graph. Let $p^{(k)}_\Gamma(v,v)$ denote the probability that the nearest neighbor random walk on $\Gamma$ starting at $v$ returns to $v$ at time~$k$. That is, $d^kp^{(k)}_\Gamma(v,v)$ is the number of walks of length $k$ in $\Gamma$ from $v$ to~$v$. We define the  \emph{spectral radius} of $\Gamma$ to be
\begin{align*}
\rho(\Gamma) \deq \sup_{v \in V}\, \limsup_{k \ra \infty} \, (p^{(k)}_\Gamma(v,v))^{1/k}.
\end{align*}
Then $\Gamma$ is non-amenable if and only if $\rho(\Gamma) < 1$ (see e.g.~\cite[$\S10$]{woess}).

Let $\pi:G \ra H$ be a surjective group homomorphism, where $G$ is some finitely generated group. We extend $\pi$ to a graph homomorphism $\pi:\Gamma_n(G) \ra \Gamma_n(H)$ given by
$$\pi(s_1, \dots, s_n) = \big(\pi(s_1), \dots, \pi(s_n)\big).$$ This is a local graph isomorphism, in other words for each $S \in \Gamma_n(G)$, the map $\pi$ induces a bijection between the edges leaving $S$ and the edges leaving~$\pi(S)$. It follows that walks in $\Gamma_n(H)$ starting at $\pi(S)$ lift uniquely to walks in $\Gamma_n(G)$ starting at~$S$. Thus, $$p^{(n)}_{\Gamma_n(G)}\big(S,S\big) \leq p^{(n)}_{\Gamma_n(H)}\big(\pi(S),\pi(S)\big),$$ and therefore $$\rho\big(\Gamma_n(G)\big) \leq \rho\big(\Gamma_n(H)\big) < 1.$$
\qed

\subsection{Proof of Lemma~\ref{lem:zz}}
The subgraph of $\Gamma_2(\zz)$ induced by $\{(a,b) \in \Gamma_2(\zz) \mid a, b > 0\}$ is a rooted binary tree, which has positive Cheeger constant. The same holds for the other three quadrants, so $\Gamma_2(\zz)$ has a subgraph $\Delta$ which is a disjoint union of four binary rooted trees. The only vertices that don't lie in $\Delta$ are $(\pm 1,0)$ and $(0, \pm 1)$, and $\Delta$ is non-amenable, so by Lemma~\ref{lem:transfer_lip}, $\Gamma_2(\zz)$ is non-amenable. By Lemma~\ref{lem:transfer_lift}, it follows that $\Gamma_2(\zz^k)$ is non-amenable for every $k \geq 1$.
\qed

\subsection{Proof of Lemma~\ref{lem:EA_zz-large}}
Let $\mathcal L$ be the class of groups which are either $\zz$-large, or finite, or infinitely generated. We want to show that every elementary amenable group belongs to~$\mathcal L$. Clearly finite groups are in~$\mathcal L$. Every finitely generated infinite abelian groups has $\zz$ as a quotient, so abelian groups are also in~$\mathcal L$. Thus by the characterization of elementary amenable groups in~\cite{chou} it is enough to show that $\mathcal L$ is closed under direct unions and extensions.

If $G$ is finite or infinitely generated, then it belongs to $\mathcal L$, so we may suppose $G$ is infinite and finitely generated. Suppose $G$ is a direct union of groups $G_i \in \mathcal L$. Since $G$ is finitely generated, $G = G_i$ for some $i$, so $G \in \mathcal L$. Now suppose $G$ is an extension of $G''$ by $G'$ where $G', G'' \in \mathcal L$. Let $\pi:G \ra G''$ be the projection with kernel~$G'$. Since $G$ is finitely generated, so is~$G''$. If $G''$ is infinite, then it must be $\zz$-large, so there is a finite index subgroup $H < G''$ which has $\zz$ is a quotient. Then $\pi^{-1}(H)$ is a finite index subgroup of $G$ which has $\zz$ as a quotient, so $G$ is $\zz$-large. On the other hand, if $G''$ is finite, then $G'$ is a finite index subgroup of $G$, so it is infinite and finitely generated. Thus it is $\zz$-large, and therefore so is~$G$.
\qed

\section{Final remarks}\label{section:final_remarks}

\subsectionn{}
Our arguments can be followed through to give explicit bounds on Cheeger constants of $\Gamma_n(G)$, and on how large $n$ must be in order for $\Gamma_n(G)$ to be non-amenable. To ease exposition, we did not track these bounds, and we did not present the arguments that would result in tight bounds. More detailed arguments with improved bounds will be presented in~\cite{thesis}.

\subsectionn{}\label{subsec:fr:gamma_schreier}
The elements of $\Gamma_n(G)$ can be identified with epimorphisms~$F_n \ra G$. The Nielsen moves then correspond to precomposition by one of the \emph{Nielsen automorphisms} $R_{ij}^{\pm 1}, L_{ij}^{\pm 1}$ of $F_n = \gen{x_1, \dots, x_n}$ given by
$$
L_{ij}^{\pm 1}(x_k) = \begin{cases}
x_k & k \neq j \\
x_i^{\pm 1}x_k & k = j,
\end{cases}
\qquad
\qquad
R_{ij}^{\pm 1}(x_k) = \begin{cases}
x_k & k \neq j \\
x_kx_i^{\pm 1} & k = j
\end{cases}
$$
with $i \neq j$ and $1 \leq i, j \leq n$. These automorphisms generate an index 2 subgroup of $\Aut(F_n)$, which we call $\Aut^+(F_n)$ (see~e.g.~\cite{whatprp, LP}). Thus, every product replacement graph is a Schreier graph of~$\Aut^+(F_n)$.

A well-known open question is whether $\Aut(F_n)$ with $n \geq 3$ has Kazhdan property $(T)$\footnote{The answer is known to be negative for $n \leq 3$~\cite{mccool,GL}.}. If $\Aut(F_n)$ has property $(T)$ for a particular value of $n$, then $\Gamma_n(G)$ is non-amenable for every infinite $n$-generated group $G$. In fact, there is a uniform lower bound on the Cheeger constants of these graphs.

\subsectionn{}\label{subsec:fr:pra}
The \emph{product replacement algorithm} is a well-known method for generating random elements of a finite group~$G$. It begins with a generating $n$-tuple $S \in \Gamma_n(G)$, and takes a random walk on $\Gamma_n(G)$, outputting a random element of the resulting generating $n$-tuple. The running time of this algorithm depends on the mixing time of the random walk on~$\Gamma_n(G)$ of~$G$. The analysis of this mixing time is also related to the question in Subsection~\ref{subsec:fr:gamma_schreier}: if $\Aut(F_n)$ has property $(T)$, then the finite product replacement graphs $\Gamma_n(G)$ form a family of expanders for any fixed $n$, and the random walk on such a graph has mixing time~$O(\log \abs{G})$. It is known that for an appropriate value of $n$, the mixing time is polynomial in $\log \abs{G}$ (see~\cite{pra_poly}). For a survey on the product replacement algorithm, see~\cite{whatprp}.

\subsectionn{}
For a finite group $G$ and a fixed number $n$, there is a lower bound on the coefficient of expansion of $\Gamma_n(G)$ in terms of the coefficients of expansion of $\Cay(G,S)$, ranging over all generating $n$-tuples $S$ (see~\cite{gamburd_pak}). Theorem~\ref{thm:uniformly_nonamenable} can be thought of as an analogue of this result for infinite groups, though the proofs differ.

\subsectionn{}\label{subsec:fr:expgrowth}
A simple consequence of Conjecture~\ref{conj:all_nonamenable} is the following.
\begin{conj}\cite{expgrowth}\label{conj:all_expgrowth}
The product replacement graph $\Gamma_n(G)$ of an infinite finitely generated group $G$ has exponential growth for sufficiently large~$n$.
\end{conj}
Some progress on this conjecture is made in~\cite{expgrowth}. Specifically, it is shown that it holds for all groups of polynomial growth, and all groups of exponential growth. It is also shown to hold for some groups of intermediate growth, including the Grigorchuk group.

\subsectionn{}
In~\cite{expgrowth} it was shown that if~$G$ has exponential growth, then $\Gamma_n(G)$ has exponential growth for every~$n \geq d(G)+1$. The proof of Theorem~\ref{thm:uniformly_nonamenable} is easily modified to prove a slight improvement of this result: if $G$ has exponential growth, then $\Gamma_n(G)$ has exponential growth for every~$n \geq d(G)$. Moreover, if $G$ has uniform exponential growth, then $\Gamma_n(G)$ also has uniform exponential growth.

\subsectionn{}
We have shown that both infinite elementary amenable groups and uniformly non-amenable groups have non-amenable product replacement graphs. A natural next step is to look at Conjecture~\ref{conj:all_nonamenable} for groups in between those two classes. This includes every group of non-uniform exponential growth:  such a group clearly cannot be uniformly non-amenable, and it has been shown that it cannot be elementary amenable either~\cite{EA_uniform_growth}. In particular, the groups of non-uniform exponential growth constructed by Wilson in \cite{wilson1, wilson2} fall between these two classes.

\subsectionn{}
Groups which are neither elementary amenable nor uniformly non-amenable belong to one of two types: amenable groups which are not elementary amenable, and non-amenable groups which are not uniformly non-amenable.

An example of the first type is the Grigorchuk group~$\gg$ (see~\cite{solved_unsolved,de_la_harpe,growthintro}). It was shown in~\cite{expgrowth} that its product replacement graphs $\Gamma_n(\gg)$ have exponential growth for $n \geq 5$, but the techniques do not appear to be strong enough to show non-amenability.

An example of the second type is the Baumslag-Solitar group $B(p,q)$ where $p$ and $q$ are relatively prime~\cite{weakly_amenable,uniform_nonamenability}. However, this group has $\zz$ as a quotient, so $\Gamma_n(B(p,q))$ is non-amenable for every $n \geq 2$. A more interesting example is the torsion group $Q$ constructed in~\cite{weakly_amenable}. Neither Theorem~\ref{thm:uniformly_nonamenable} nor Theorem~\ref{thm:zz_large} is enough to show that $\Gamma_n(Q)$ is non-amenable for some $n$.

\subsectionn{}
Another example of interest is Thompson's group~$F$ (see~\cite{thompsons_group}). Whether $F$ is amenable is a well-known open problem, but it is known that it is not elementary amenable. Thompson's group $F$ has~$\zz^2$ as a quotient, and therefore $\Gamma_n(F)$ is non-amenable for every~$n \geq 2$. However, the related groups $T$ and $V$ are both simple groups. Thus, they cannot be $\zz$-large, and Theorem~\ref{thm:zz_large} does not apply. Note that $T$ and $V$ both have exponential growth, so by the results in \cite{expgrowth} they satisfy Conjecture~\ref{conj:all_expgrowth}.

\section{Acknowledgements}

I am grateful to my graduate advisor Igor Pak for suggesting the problem of non-amenability of infinite product replacement graphs. I would like to thank Yehuda Shalom for pointing out what is known about uniform non-amenability.


\begin{thebibliography}{Dun70}

\bibitem[A+]{uniform_nonamenability}
G.~N.~Arzhantseva, J.~Burillo, M.~Lustig, L.~Reeves, H.~Short and E.~Ventura,
Uniform non-amenability,
\emph{Adv. Math.}~\textbf{197} (2005), 499--522.

\bibitem[BL]{BL}
H.~B\"{a}\"{a}rnhielm and C.~R.~Leedham-Green,
The product replacement prospector,
\emph{J.~Symbolic Comput.}~\textbf{47} (2012), 64--75.

\bibitem[BLM]{BLM}
A.~V.~Borovik, A.~Lubotzky and A.~G.~Myasnikov,
The finitary Andrews-Curtis conjecture,
in \emph{Progr. Math.}~\textbf{248},
Birkh\"{a}user, Basel, 2005, 15--30.

\bibitem[BG]{linear_independence}
E.~Breuillard and T.~Gelander,
Uniform independence in linear groups,
\emph{Invent. Math.}~\textbf{173} (2008), 225--263.

\bibitem[CFP]{thompsons_group}
J.~W.~Cannon, W.~J.~Floyd and W.~R.~Parry,
Introductory notes on Richard Thompson's groups. 
\emph{Enseign. Math.}~\textbf{42} (1996), 215--256.

\bibitem[Chou]{chou}
C.~Chou,
Elementary amenable groups,
\emph{Illinois J. Math.}~\textbf{24} (1980), 396--407.

\bibitem[dlH]{de_la_harpe}
P.~de~la Harpe, \emph{Topics in Geometric Group Theory}, University of Chicago
  Press, Chicago, 2000.

\bibitem[EP]{largeness}
M.~Edjvet and S.~J.~Pride,
The concept of ``largeness'' in group theory II,
in \emph{Lecture Notes in Math.,}~\textbf{1098}, Springer, Berlin, 1984, 29--54.

\bibitem[Eva]{evans}
M.J.~Evans,
Nielsen equivalence classes and stability graphs of finitely generated groups, in \emph{Ischia group theory 2006}, World Sci. Publ., Hackensack, 2007, 103--119.

\bibitem[GP]{gamburd_pak}
A.~Gamburd and I.~Pak,
Expansion of product replacement graphs,
\emph{Combinatorica}~\textbf{26} (2006), 411--429.

\bibitem[Gri]{solved_unsolved}
R.~I.~Grigorchuk, Solved and unsolved problems around one group, in
  \emph{Infinite Groups: Geometric, Combinatorial and Dynamical Aspects},
  Birkh\"auser, Basel, 2005, 117--218.

\bibitem[GP]{growthintro}
R.~I.~Grigorchuk and I.~Pak, Groups of intermediate growth, an introduction,
  \emph{L'Ens. Math.}~\textbf{54} (2008), 251--272.

\bibitem[Gro]{gromov}
M.~Gromov, Groups of polynomial growth and expanding maps,
\emph{IHES~Publ.~Math.}~\textbf{53} (1981), 53--78.

\bibitem[GL]{GL}
F.~Grunewald and A.~Lubotzky,
Linear representations of the automorphism group of a free group,
\emph{Geom. Funct. Anal.}~\textbf{18} (2009), 1564--1608.

\bibitem[LP]{LP}
A.~Lubotzky and I.~Pak,
The product replacement algorithm and Kazhdan's property~$(T)$,
\emph{J.~AMS}~\textbf{14} (2001), 347--363.

\bibitem[LZ]{property_tau}
A.~Lubotzky and A.~\.{Z}uk, \emph{On property $(\tau)$}, monograph in preparation.

\bibitem[Mal]{thesis}
A.~Malyshev, \emph{Combinatorics of finitely generated groups}, Ph.D.~thesis,
  UCLA, in preparation.

\bibitem[MP]{expgrowth}
A.~Malyshev and I.~Pak, Growth in product replacement graphs, {\tt arXiv:}{\tt 1304.5320}.

\bibitem[Mc]{mccool}
J.~McCool,
A faithful polynomial representation of $\Out F_3$. 
\emph{Math. Proc. Cambridge Philos. Soc.}~\textbf{106} (1989), 207--213. 

\bibitem[Met]{Met}
W.~Metzler, On the Andrews-Curtis conjecture and related problems, in
\emph{Contemp. Math.}~\textbf{44}, AMS, Providence, RI, 1985, 35--50.

\bibitem[Myr]{myropolska}
A.~Myropolska,
Andrews--Curtis and Nielsen equivalence relations on some infinite groups,
{\tt arXiv:}{\tt 1304.2668}.

\bibitem[O1]{weakly_amenable}
D.~V.~Osin, Weakly amenable groups.
\emph{Contemp. Math.}~\textbf{298} (2002), 105--113.

\bibitem[O2]{EA_uniform_growth}
D.~V.~Osin,
Algebraic entropy of elementary amenable groups,
\emph{Geom. Dedicata}~\textbf{107} (2004), 133--151. 

\bibitem[P1]{whatprp}
I.~Pak, What do we know about the product replacement algorithm?, in
  \emph{Groups and Computation III}, de Gruyter, Berlin, 2001, 301--347.

\bibitem[P2]{pra_poly}
I.~Pak, The product replacement algorithm is polynomial,
in \emph{Proc. FOCS 2000}, IEEE Comput. Soc. Press,
Los Alamitos, CA, 2000, 476--485.

\bibitem[W1]{wilson1}
J.~S.~Wilson,
On exponential growth and uniformly exponential growth for groups, 
\emph{Invent. Math.}~\textbf{155} (2004), 287--303. 

\bibitem[W2]{wilson2}
J.~S.~Wilson,
Further groups that do not have uniformly exponential growth,
\emph{J. Algebra}~\textbf{279} (2004), 292--301. 

\bibitem[Woe]{woess}
W.~Woess, \emph{Random walks on infinite graphs and groups}, Cambridge U.~Press, Cambridge, 2000. 

\end{thebibliography}
\end{document}